\documentclass[reqno,12pt]{amsart}

\usepackage{epsf}
\usepackage{graphics}
\usepackage{amssymb}
\usepackage{amsmath}

\date{}

\theoremstyle{plain}
\newtheorem{theorem}{Theorem}
\newtheorem{corollary}{Corollary}
\newtheorem{lemma}{Lemma}

\theoremstyle{definition}

\theoremstyle{remark}

\newtheorem*{remark}{Remark}

\def\C{{\mathbb C}}

\def\R{{\mathbb R}}
\def\Z{{\mathbb Z}}

\title{Knotted holomorphic discs in $\C^2$} 

\author{S.~Baader, F.~Kutzschebauch, E.~F.~Wold}

\begin{document}

\begin{abstract} We construct knotted proper holomorphic embeddings of the unit disc in $\C^2$.
\end{abstract}

\maketitle

\section{Introduction}

Every classical knot type can be represented by a polynomial embedding of 
$\R$ in $\R^3$~\cite{S}. In particular, there exist topologically distinct polynomial 
embeddings of $\R$ in $\R^3$. Crossing these with another coordinate, we obtain 
topologically distinct polynomial embeddings of $\R^2$ in $\R^4$. In contrast, all 
polynomial embeddings of $\C$ in $\C^2$ are topologically equivalent, in fact even 
connected by a polynomial automorphism of $\C^2$. 
The first complete algebraic proof of this fact is due to Abhyankar and Moh~\cite{AM}; 
later a purely knot theoretical proof was found by Rudolph~\cite{R1}. It is an open 
question whether proper holomorphic embeddings of the complex plane $\C$ 
or of the unit disk $\triangle\subset\mathbb C$
in $\C^2$ can be topologically knotted 
(see Problem 1.102 (A/B) of Kirby's list~\cite{K}).
In this note we construct knotted proper holomorphic embeddings of the unit disc in~$\C^2$. 

\begin{theorem} \label{mainthm}
There exist topologically knotted proper holomorphic embeddings of the unit disc in $\C^2$.
\end{theorem}

Since the image of a holomorphically embedded disc in $\C^2$ is a minimal surface, we obtain the following corollary, which solves Problem 1.102 (C) of Kirby's list.

\begin{corollary} 
There exists a proper embedding $f: \R^2 \rightarrow \R^4$ whose image is a topologically knotted complete minimal surface.
\end{corollary}

Our construction is based on the existence of locally well-behaved Fatou-Bieberbach domains 
and on the existence of knotted holomorphic discs in the 4-ball. A Fatou-Bieberbach domain 
in $\C^2$ is an open subset of $\C^2$ which is biholomorphically equivalent to $\C^2$.
Fatou-Bieberbach domains tend to have wild shapes; the first Fatou-Bieber\-bach domain 
with smooth boundary was constructed by Stens\o nes~\cite{St}. In~\cite{G}, 
Globevnik constructed Fatou-Bieberbach domains $\Omega \subset \C^2$ whose 
intersections with $\C \times \Delta$ are arbitrarily small $\mathcal{C}^1$-perturbations 
of $\Delta \times \Delta$. Here $\Delta$ stands for the open unit disc in $\C$. 
These domains are well-adapted for constructing knotted holomorphic discs in $\C^2$.
Indeed, all we need is a knotted proper holomorphic embedding of the closed unit disc 
$\varphi: \overline{\Delta} \rightarrow \overline{\Delta} \times \overline{\Delta}$ that maps $\partial \overline{\Delta}$ to $\partial \overline{\Delta} \times \Delta$. Composing the restriction of this embedding $\varphi |_\Delta: \Delta \rightarrow \Delta \times \Delta \subset \Omega$ with a biholomorphism $h: \Omega \rightarrow \C^2$ yields a knotted proper holomorphic embedding of $\Delta$ in~$\C^2$. The existence of knotted proper holomorphic embeddings $\varphi: \overline{\Delta} \rightarrow \overline{\Delta} \times \overline{\Delta}$ is easily established by using the theory of complex algebraic curves in $\C^2$. 

We describe the relevant features of the Fatou-Bieberbach domains needed for our purpose
in Section~2. The proof of Theorem~\ref{mainthm} is completed in Section~3, 
where we give an explicit example of a knotted holomorphic disc in $\overline{\Delta} 
\times \overline{\Delta}$.

\section{Fatou-Bieberbach domains}

Fatou-Bieberbach domains in $\C^2$ are usually described by certain infinite processes, for example as domains of convergence of maps defined by sequences of automorphisms. This technique was used by Globevnik to construct Fatou-Bieberbach domains with controlled shape inside $\C \times \Delta$.

\begin{theorem}[Globevnik \cite{G}] \label{globevnik}
Let $Q \subset \C$ be a bounded open set with boundary of class 
$\mathcal{C}^1$ whose complement is connected. Let $0<R<\infty$ be such that 
$\overline{Q} \subset R \Delta$. There are a domain $\Omega \subset \C^2$ and a 
volume-preserving biholomorphic map from $\Omega$ onto $\C^2$ such that

\smallskip
\begin{enumerate}
\item[(i)] $\Omega \subset \{(z,w) \in \C^2:|z|<\mathrm{max} \{R,|w|\}\}$, \\
\item[(ii)] $\Omega \cap R(\Delta \times \Delta)$ is an arbitrarily small $\mathcal{C}^1$-perturbation of $Q \times R \Delta$. \\
\end{enumerate}
\end{theorem}

The assumptions of Theorem~\ref{globevnik} are verified when $Q$ is an open disc with smooth boundary (and $R>0$ large enough). This is the version we will need. As mentioned in the introduction, we will insert knotted discs in $\Delta \times \Delta$ with boundary in $\partial \overline{\Delta} \times \Delta$. It is a priori not clear whether such discs stay knotted in the larger domain $\Omega$. The following lemma allows us to control knottedness on the level of the fundamental group.

\begin{lemma} \label{lemma1}
Let $\Omega \subset \C^2$ be an open domain homeomorphic to $\C^2$ with $
\Omega \cap \C \times \Delta=\Delta \times \Delta$, and let $X \subset \Delta 
\times \Delta$ be a subset with $\overline{X} \subset \overline{\Delta} \times 
\Delta $. The inclusion $i: (\Delta \times \Delta) \setminus X \rightarrow \Omega 
\setminus X$ induces an injective map 
$i_*: \pi_1 ((\Delta \times \Delta) \setminus X) \rightarrow \pi_1 (\Omega 
\setminus X)$ (suppressing base points).
\end{lemma}

\begin{remark} At this point it does not matter whether $\Omega \cap \C \times \Delta$ 
is precisely $\Delta \times \Delta$ or a small $\mathcal{C}^1$-perturbation of it. 
Further, the corresponding statement stays true if the intersection $\Omega \cap (\C \times 
\Delta)$ is a small $\mathcal{C}^1$-perturbation of $D \times \Delta$, where 
$D \subset \C$ is any embedded disc with smooth boundary.
\end{remark}

\begin{proof} Choose $\epsilon>0$ such that $X \subset \Delta \times (1-\epsilon) \Delta$. Setting 
$$U:=(\Delta \times \Delta) \setminus X$$ 
and 
$$V:=\Omega \setminus (\Delta \times \overline{(1-\epsilon) \Delta}),$$ 
we have $\Omega \setminus X=U \cup V$. The theorem of Seifert-van Kampen tells us 
that $\pi_1 (\Omega \setminus X)$ is the free product of $\pi_1 (U)$ and $\pi_1 (V)$
 amalgamated over $\pi_1 (U \cap V)$. The latter is isomorphic to $\Z$, since 
$$U \cap V=\Delta \times (\Delta \setminus \overline{(1-\epsilon) \Delta})$$ 
is homotopy equivalent to a circle. Let $\gamma: S^1 \rightarrow U \cap V$ be a loop 
that generates $\pi_1 (U \cap V)$. If the induced map $i_*: \pi_1 ((\Delta \times \Delta)
\setminus X) \rightarrow \pi_1 (\Omega \setminus X)$ were not injective, then a 
certain non-zero multiple of $[\gamma] \in \pi_1 (V)$ would have to vanish. 
This is impossible since the inclusion $j: V \rightarrow \C \times \C^*$ maps $\gamma$ 
onto a generator of $\pi_1 (\C \times \C^*) \cong \Z$.
\end{proof}

\section{Complex plane curves}

A complex plane curve is the zero level $V_f=\{(z,w) \in \C^2: f(z,w)=0\} \subset \C^2$ of a non-constant polynomial $f(z,w) \in \C[z,w]$. Complex plane curves form a rich source of examples in algebraic geometry and topology. For example, the intersection of a complex plane curve $V_f$ with the boundary of a small ball centred at an isolated singularity of $f$ forms a link which is often called algebraic. The class of algebraic links has been generalized by Rudolph to the larger class of quasipositive links~\cite{R2}. Using the theory of quasipositive links, it is easy to construct knotted proper holomorphic embedded discs in $\Delta \times \Delta$. In the following, we will give a short description of Rudolph's theory; more details are contained in~\cite{R2} and~\cite{R3}.

Let $f(z,w)=f_0(z)w^n+f_1(z)w^{n-1}+\ldots+f_n(z) \in \C[z,w]$ be a non-constant polynomial, $f_0(z) \ne 0$. Under some generic conditions on $f$, the set $B$ of complex numbers $z$ such that the equation $f(z,w)=0$ has strictly less than $n$ solutions $w$, is finite. Let $\gamma \subset \C \setminus B$ be a smooth simple closed curve. The intersection $L=V_f \cap \gamma \times \C$ is a smooth closed $1$-dimensional manifold, i.e. a link, in the solid torus $\gamma \times \C$. More precisely, the link $L$ is an $n$-stranded braid in $\gamma \times \C$, which becomes a link in $S^3$ via a standard embedding of $\gamma \times \C$ in $S^3$. Links that arise in this way are called quasipositive.

If $D \subset \C$ is the closed disc bounded by $\gamma$, then the intersection $X=V_f \cap D \times \C$ is a piece of complex plane curve bounded by the link $L$. Choosing the polynomial $f$ and the curve $\gamma$ appropriately, it is possible to arrange $X$ to be an embedded disc with a non-trivial quasipositive knot as boundary (see~\cite{R2}, Example 3.2). Assuming that $X$ is compact, there exists $R>0$, such that $X \subset D \times R \Delta$. In order to establish Theorem~\ref{mainthm}, it remains to find an example where $X$ is knotted in $D \times R \Delta$, more precisely
$$\pi_1(D \times R \Delta \setminus X) \ncong \Z.$$
In view of Lemma~\ref{lemma1}, this will then give rise to a knotted proper holomorphic embedding of the unit disc in $\Omega$, hence in $\C^2$.

\begin{remark} Here we choose $\Omega \subset \C^2$ to be a Fatou-Bieberbach domain 
whose intersection $\Omega \cap \C \times R \Delta$ is a small $\mathcal{C}^1$-
perturbation of $D \times R \Delta$. If this perturbation is small enough, then the 
intersection $\widetilde{X}=V_f \cap \C \times R \Delta \cap \Omega$ is still a proper
holomorphic embedded disc, as knotted as $X$ in $D \times R \Delta$.
\end{remark}

We conclude this section with a concrete example, in fact Rudolph's Example~3.2 
of~\cite{R2}. Let
$$f(z,w)=w^3-3w+2z^4.$$
The complex plane curve $V_f$ is easily seen to be non-singular. The equation $f(z_0,w)=0$ 
fails to have three distinct solutions $w$, if and only if the two equations $f(z_0,w)=0$ and 
$\frac{\partial}{\partial w}f(z_0,w)=3w^2-3=0$ have a simultaneous solution. This 
happens precisely when $z_0^8=1$. Thus the set $B$ consists of the $8$th roots of unity.
For $z \in \C \setminus B$, the equation $f(z,w)=0$ has three distinct solutions $w_1, 
w_2, w_3 \in \C$, which we index by increasing real parts. Further, it will be convenient to 
determine the set $B_+$ of complex numbers $z$ such that the equation $f(z,w)=0$ has two 
distinct solutions $w$ with coinciding real parts. In our example, $B_+$ consists of $8$ rays 
emanating from the points of $B$, as shown in Figure~1. These rays carry labels $1$ or $2$, 
depending on whether the real parts of $w_1$ and $w_2$ or those of $w_2$ and $w_3$ 
coincide. There is a way of orienting the rays that corresponds to choosing positive standard 
generators of the braid group. Figure~1 also indicates a curve $\gamma$ that bounds a disc~
$D$. We claim that the piece of complex plane curve $X=V_f \cap D \times \C$ is a disc
with $\pi_1((D \times \C) \setminus X) \ncong \Z$.

\begin{figure}[ht]
\scalebox{1}{\raisebox{-0pt}{$\vcenter{\hbox{\epsffile{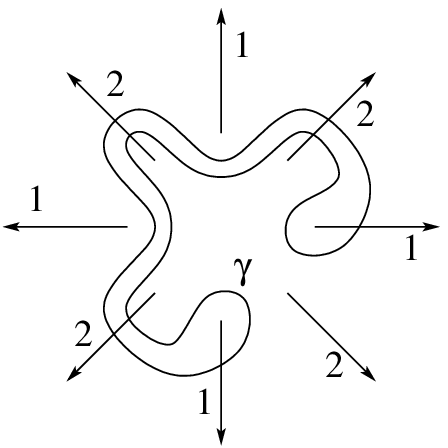}}}$}} 
\caption{}
\end{figure}

First, we observe that if we `cut off' the two ends of the disc $D$ of Figure~1, we obtain 
another disc $D'$ disjoint from $B$. The intersection $V_f \cap D' \times \C$ is therefore a 
disjoint union of three discs. Adding the ends to $D'$ again gives rise to two identifications 
along the boundaries of these discs. It is easy to see that these identifications result in a single 
disc whose boundary $L=V_f \cap (\gamma \times \C)$ is knotted (actually $L$ is the 
quasipositive ribbon knot $8_{20}$). However, the mere fact that $L$ is knotted does not imply 
that the fundamental group of $(D \times \C) \setminus X$ 
is not isomorphic to $\Z$ (see the last paragraph of~\cite{FT}). We will describe 
$\pi_1((D \times \C) \setminus X)$ explicitly, in terms of generators and relations. Hereby 
we will use Orevkov's method for presenting the fundamental group of the complement of 
complex plane curves in $\C^2$~\cite{O}. To every connected component of 
$\C \setminus B_+$ (in our case $D \setminus B_+$), we assign $n$ (in our case $3$) 
generators corresponding to meridians of the $n$ discs lying over that component. Every 
oriented edge of the graph $B_+$ gives rise to $n$ relations between the generators of the 
two (not necessarily distinct) components adjacent to that edge. All these relations are of 
Wirtinger-type (see~\cite{O}, Lemma~3.1).

Let us denote the four connected component of $D \setminus B_+$ by $U_1$, $U_2$, $U_3$, $U_4$, starting at the bottom left and going clockwise around $D$. To every region $U_j$ correspond three generators $\alpha_{ij}=\alpha_i(U_j)$, $1 \leq i \leq3$. The relations among these generators read as follows:

\smallskip
1. edge: $\alpha_{11}=\alpha_{21}$ 

\smallskip
2. edge: $\alpha_{11}=\alpha_{12}$, $\alpha_{31}=\alpha_{22}$, $\alpha_{21}=\alpha_{31} \alpha_{32} \alpha_{31}^{-1}$ 

\smallskip
3. edge: $\alpha_{12}=\alpha_{13}$, $\alpha_{32}=\alpha_{23}$, $\alpha_{22}=\alpha_{32} \alpha_{33} \alpha_{32}^{-1}$ 

\smallskip
4. edge: $\alpha_{13}=\alpha_{14}$, $\alpha_{33}=\alpha_{24}$, $\alpha_{23}=\alpha_{33} \alpha_{34} \alpha_{33}^{-1}$ 

\smallskip
5. edge: $\alpha_{14}=\alpha_{24}$ 

\smallskip
\noindent 
From this it is easy to verify that the assignment $\alpha_{11} \mapsto (23)$, 
$\alpha_{31} \mapsto (12)$ defines a surjective homomorphism 
$\varphi: \pi_1((D \times \C) \setminus X) \rightarrow S_3$, whence 
$\pi_1((D \times \C) \setminus X)$ is not a cyclic group.

\bigskip
\address{ETH Z\"urich, R\"amistrasse~101, CH-8092 Z\"urich, Switzerland}

\email{\texttt{sebastian.baader@math.ethz.ch}}

\bigskip
\address{Universit\"at Bern, Sidlerstrasse~5, CH-3012 Bern, Switzerland}

\email{\texttt{Frank.Kutzschebauch@math.unibe.ch}}

\bigskip
\address{Universitetet~I Oslo, Postboks 1053 Blindern, 0316 Oslo, Norway}

\email{\texttt{erlendfw@math.uio.no}}

\end{document}